\UseRawInputEncoding

\documentclass[12pt]
{amsart}
\usepackage[russian]{babel}
\usepackage[cp1251]{inputenc}
\usepackage{amsfonts,amsmath,amssymb,amscd,amsthm,url,graphicx}
\long\def\comment#1\endcomment{}

\usepackage[X2,T2A]{fontenc}
\DeclareTextSymbolDefault{\cyryat}{X2}
\DeclareTextSymbolDefault{\CYRYAT}{X2}

\theoremstyle{theorem}
\newtheorem{theorem}{Теорема}
\theoremstyle{definition}
\newtheorem{blud}[theorem]{Заблуждение}

\begin{document}

%см. (8) в обосновании некомпетентности

%<div><p></p></div>
%Research projects and Moscow Mathematical Conference for high school students
%A. Zaslavskiy and A. Skopenkov
%This paper shares some experience in advanced mathematical education. 
%We show how a high school student can be naturally and gradually introduced to basic steps of scientific research: developing intuition by finding and correcting mistakes through discussions and writing a paper,  (transparent) anonymous  peer review, recognition and award.
%We show that most of this can be done in research projects not aiming at scientific novelty.
%%which do not claim any scientific novelty (or even are explicitly known).
%We share the experience (both principles and examples) of the Moscow Mathematical Conference of High School Students.  

%impartiality principle

%\newpage
\title{Исследовательские задачи и Московская математическая конференция школьников}

\author{А.А. Заславский и А.Б. Скопенков}

%Эта версия представлена в <<Математическое Просвещение>> в январе 2026 г.
%Предыдущие версии предварительно обсуждались с августа 2024 г. (с ответственным секретарем редколлегии и другими деятелями математического образования) и   
%представлялись в МП в феврале 2025 г. и в ноябре 2025.  

\thanks{Благодарим Д.В. Прокопенко, А.И. Сгибнева, М.Б. Скопенкова, Б.Р. Френкина и Г.Р. Челнокова за полезные замечания и обсуждения.
О нарушениях принципа непредвзятости и профессиональной этики при рецензировании редколлегией <<Математического Просвещения>> предыдущих версий данной статьи  см. \cite{ZS26}. 
\newline
А.А. Заславский: ЦЭМИ РАН, МЭИ; \texttt{alzasl@yandex.ru}.
\newline
А.Б. Скопенков: \url{https://users.mccme.ru/skopenko/};
Московский физико-технический институт и Независимый Московский университет; \texttt{skopenko@mccme.ru}.}

\date{}

\maketitle

\tableofcontents

%\abstract\endabstract

\section*{Введение}\label{s:intr}

Российское (и советское) образование знаменито <<системой кружков и олимпиад>>, включающей также спецкурсы, <<спецматематику>>, летние школы, подготовку к олимпиадам и т.д.
В этом тексте мы поделимся опытом подкрепления и развития этой системы в рамках Московской математической конференции школьников.
Далее мы называем ее просто Конференцией, хотя она включает также конкурс и награждение работ.
Мы покажем, как преодолеть широко известные недостатки многих конференций-конкурсов школьников: 
работы не самостоятельные, результаты недостаточно проверяются, доклады излагаюстся непонятно.
Такие недостатки подталкивают школьников к \emph{имитации} исследовательской деятельности.
Мы надеемся, что этот текст будет интересен участникам <<системы кружков и олимпиад>> (школьникам, учителям и исследователям), организаторам конференций школьников, а также людям, интересующимся математическим образованием как частью культуры.

Мы проводим Конференцию, чтобы приобщить школьников к решению \textbf{исследовательских задач}.
К ним мы относим задачи, в которых ответ или метод решения неизвестен ученику (но не обязательно нерешенные научные задачи).
Важное отличие исследовательской деятельности от  учебной --- необходимость более серьезно проверять результаты (\S\ref{zaver}\footnote{\label{f:editors}Читатель может читать статью подряд --- и тогда пропускать ссылки.
Ссылки приведены, чтобы помочь читателю поскорее добраться до интересующего его материала, если он захочет сделать это вместо чтения подряд.}). 
%[Эта сноска добавлена в соответствии с замечанием редколлегии Мат. Просвещения <<\emph{...сложная система отсылок в разные разделы, читать тяжело...}>>. 
%Мы согласны с ответственным секретарем редколлегии, что информация в сноске само собой разумеющаяся и сноску желательно удалить. 
%Мы будем рады это сделать, если редколлегия подтвердит, что соответствующее замечание снято (даже при отсутствии сноски) и не будет далее повторяться.]  
%Ответственный секретарь сообщил, что редколлегия не способна дать такое подтверждение в явной форме. 
%Тогда его можно дать в форме решения о принятии статьи с данными конкретными сокращениями.
%?В чем заключается подкрепление и развитие <<системы кружков и олимпиад>>?
\emph{Мы приобщаем учеников к идеям и развиваем их интуицию.}
Мы включаем в этот процесс проверку записью --- более серьезную, чем на олимпиадах. 
%{\bf ?} некоторых {\bf ?} многих [arxiv современных] 
%Мы не исключаем из
Здесь под интуицией мы понимаем именно способность к озарению, а не само озарение.
Грубо говоря, если я не только продумал решения десяти задач на интуитивном уровне, но и записал их, то моя интуиция при решении одиннадцатой задачи уже выше.

Развитие интуиции реализуется профессиональным общением с консультантом и рецензентом (\S\ref{s:rest}).
%приобщение школьников к \emph{учебно-исследовательской} работе (\S\ref{why}, \S\ref{s:kaki}).
%А для школьников, достигших соответствующего уровня, ---
\emph{Консультации --- важнейшая форма работы,} но она наименее формализована.
%О ней труднее всего написать.
В результате консультаций не только обнаруживаются и исправляются ошибки, но рождаются упрощения (см. пример работы К. Зюбина в \S\ref{s:examp})
%формулировок и доказательств,
и новые результаты (см. пример работы В. Зеленина в \S\ref{s:examp}).

Это 
%[развитие] 
происходит и в процессе записи, и в неформальных обсуждениях
(возможность которых основана на более раннем серьезном  освоении базового материала).
Основы этой работы описаны в 
\S\ref{s:howto}, \S\ref{zaver} и в конце \S\ref{s:rest}.

%https://old.mccme.ru//circles//oim/mmks/works2022/suvorov1.pdf
%https://old.mccme.ru//circles//oim/mmks/works2022/suvorov4.pdf
%Или: 2023, Хоменко Анастасия и Кеелус Милена

На продвинутом уровне этот процесс формирует реалистичное представление о \emph{научно-ис\-сле\-до\-ва\-тельской} работе.
Самостоятельное решение учебных задач (включающее проверку путем записи) создает ту базу, из которой только и  
%(по нашему мнению, и только из которой) 
рождаются впоследствии открытия (не обязательно в математике!).
Например, Г. Перельман начинал с решения учебных задач на кружке (не важно,что внешней мотивировкой кружка, в котором он занимался, были письменные олимпиады, а не конференции школьников).
Надеемся, что участие в Конференции поможет школьникам создать базу для профессионального роста.

%Хорошо, что Вы не ставите задачу быстро сделать научную работу.
%А видите ее в том, чтобы потренироваться грамотно записывать математические тексты.

Научно-педагогическому сообществу полезно и интересно получить представление о <<научно-проектно-конференционной>> деятельности школьников в целом.
Пока что это невозможно, ибо результаты этой деятельности недоступны: из многих известных нам конференций школьников полные тексты работ выкладываются в интернет только на ММКШ.
См. подробнее <<Прозрачность>> в \S\ref{s:rest}. 

Подготовка работы на Конференцию --- хорошая тренировка к олимпиадам.
Действительно, подобно решению задачи на олимпиаде, существенная часть такой подготовки --- придумывание и написание четких формулировок и завершенных доказательств (\S\ref{zaver}).
В процессе такой подготовка ученик развивается в более <<взрослом>> направлении, чем на олимпиадах.
Это помогает (а не мешает) подготовке к самим олимпиадам (парадокс? читайте \cite{Sk}, \cite[\S26]{ZSS}).
В частности, доведение некоторых доказательств до \emph{завершенных}, уровень надежности которых превышает уровень <<полного балла на олимпиаде>>, так же полезно, как тренировка на песке перед футбольным матчем на траве.

Многие из приводимых соображений не оригинальны.
Однако опыт авторов говорит о необходимости еще раз высказать эти соображения.\footnote{{\it Toutes choses sont dites d\'ej\`a; mais comme  personne n'\'ecoute, il faut toujours recommencer. A. Gides.} 
\newline
Все, что должно быть сказано, уже сказано. Но поскольку никто не слушал, нужно постоянно начинать заново. Андре Жид.} 
% (вольный перевод авторов)

%Ach, alle Weisheit ist so einfach, ist schon so lange, schon so genau und unzweifelhaft ausgesprochen und formuliert worden! Warum geh\"ort sie uns nur zuzeiten, nur an den guten Tagen, warum nicht immer? H. Hesse.
%Voir sur https://citations.ouest-france.fr/citation-andre-gide/choses-sont-dites-deja-comme-30118.html

См. также расширенную обновляемую версию \cite{P} этой статьи и \cite{MM, Sg13, Sg, Sk08, SS}.
О других формах <<проектно-конференционной>> деятельности школьников см., например, \cite{LK, RPT} и перечисление в \cite[конец параграфа 2 (<<позитивный опыт>>)]{Be14}.

\section{Зачем участвовать в Конференции и как обычно возникают работы?}\label{why}

Изучение математики путем решения и разбора задач --- важнейшая черта <<системы кружков и олимпиад>>, см. подробнее \cite[п. 1.2]{ZSS}.
Мы называем такую деятельность \textbf{учебно-исследовательской}.
Прилагательным <<исследовательская>> мы подчеркиваем

$\bullet$ самостоятельность и творческий характер (изучения материала путем) решения задач (пусть даже несложных);

$\bullet$ естественность постепенного перехода к научно-исследовательской работе.

При таком изучении часто возникают решения, которые интересно показать более широкому кругу школьников и учителей.
Иногда при этом естественно появляется элемент новизны или даже новые результаты.
Но это совсем не обязательно: такие решения могут быть и самостоятельным переоткрытием известной теоремы.
Их разумно и полезно подавать на Конференцию.
% (можно перед этим направлять на консультацию, см. \S\ref{s:cons}).
В частности, удачно подавать на Конференцию решение задачи, предложенной на Летней конференции Турнира городов \cite{LK} (далее --- ЛКТГ), в Задачнике <<Кванта>> или на самой Конференции \cite{M}.
%Особенно если это решение отличается от предложенного жюри/редакцией.
См. также \cite[Как найти задачу для исследования, стр. 25-26]{Sk08}.

Обучение проверке идей путем записи решений --- процесс не только интересный и полезный, но длительный и трудоемкий.
Принятие работы на Конференцию (и, если автору удобно, само выступление) вносит в него элемент праздника.
%дополнительно его стимулирует, что и является одной из функций Конференции.
%(Даже если автор не смог выступить, само принятие работы его мотивирует.)

Перед участниками выступают и школьники, и математики (не только московские!).
Школьники рассказывают о своих результатах (\S\ref{s:examp}, \cite{N}), а математики представляют задачи для исследования \cite{M}.
Заседания проходят живо: высказываются интересные комментарии и задаются вопросы.
Доклады понятны школьникам и учителям; даже в сложных докладах приводятся доступные им яркие результаты и идеи доказательств 
%доступные широкой аудитории.
См., например, видеозаписи некоторых докладов \cite{N}.
%Результаты докладов либо надёжно проверены , либо явно названы гипотезами.
Участие в Конференции помогло многим авторам довести работы до публикации в <<Кванте>>, в <<Математическом Просвещении>> и даже в научных журналах \cite{Ko17, KP, Mo16, Ni08, Ru10}.
%Вот примеры работ: \S\ref{s:examp} (характерные) и \cite{N}.
%Что означает фраза (все работы второго очного тура).
%См. также \cite{R}.

%Задачи для исследования <a href=http://www.mccme.ru/circles/oim/mmks/notesm.htm> предлагаются и на Конференции</a>.
%Однако совсем не обязательно, чтобы в дальнейшем школьники, пришедшие на Конференцию, занимались решением именно этих задач.

За счет чего это получается?
%О работе школьника перед подачей на Конференцию см. \S\ref{why}.
Поданные тексты тщательно рецензируются (\S\ref{s:rest}).
Сначала рецензируются краткие версии, содержащие только четкую формулировку одного результата с необходимыми определениями (см. подробнее \S\ref{s:howto} и \cite[\S0]{P}).
Авторы учитывают замечания рецензентов, по желанию пользуясь помощью консультантов (которые обычно не становятся соруководителями).

%Иногда такая доработка превращается в содержательное переделывание работы.
%Когда доработка завершена, текст принимается на Конференцию и выкладывается на её сайт.
%Такой серьёзный труд интересен увлечённому математикой школьнику и вознаграждается вниманием к его результату (а также премией).

\section{Какие работы принимаются на Конференцию?}\label{s:kaki}

Мы принимаем \emph{самостоятельные содержательные} (в т.ч. простые и короткие) работы (=тексты)
школьников и первокурсников одного из следующих типов (\emph{номинаций}):

$\bullet$ научно-исследовательские работы (есть четкая формулировка и завершенное доказательство, см. \S\ref{zaver}; работа представлена в Архив \url{http://arxiv.org} или в \url{https://hal.science}, cм. подробнее начало \S\ref{s:nomi}).
%появление на математическом сайте\footnote{\label{f:arxiv} Более конкретно, нужно представить работу в один из архивов \url{http://arxiv.org} или \url{https://hal.science}.

%Если появление представленной работы в Архиве задерживается (или работа отклонена модераторами Архива), то по просьбе автора и по решению программного комитета требование появления работы в Архиве до 15.11
%может быть заменено на требование появления работы на другом  до 15.11.

$\bullet$ учебно-исследовательские работы (есть четкая формулировка и завершенное доказательство);

$\bullet$ исследовательские разработки (есть четкая формулировка гипотезы);

$\bullet$ наглядные / экспериментальные материалы (есть понятные рисунок или видео / ссылка на архив типа github, на котором лежит программа и математически грамотное описание ее ввода и вывода).

%[того, что она принимает на вход и что выводит]).
%, доступной для тестирования, программа проходит тесты,
%(на сайт Конференции выкладывается не сама программа, а описание и ссылка)

Работа не должна содержать ни одного заведомо неверного утверждения, а также необоснованной претензии на новизну (результатов или доказательств) 
или на завершенность доказательств. 
(Если это не выполнено, то автору предоставляется возможность подать новую версию, в которой таковых не будет.)

Во все номинации, кроме первой, принимаются результаты как заведомо не претендующие на новизну (например, решение задачи с кружка), так и потенциально имеющие элемент новизны (которая пока не проверена, см. начало \S\ref{s:nomi}).

Эти критерии (а также приводимые ниже и в \cite[\S2]{P} пояснения к ним) не очень формальны.
Решение о принятии / отклонении работы по каждой номинации принимает программный комитет на основании рецензий, см. \S\ref{s:rest}.
Он же распределяет время на доклады, а доклады --- на аудиторные и стендовые.
Для последних решений важны степень содержательности и интересности работы (а не только завершенность доказательств).

\textbf{Примеры работ,} принятых в каждую из номинаций, см. в \S\ref{s:examp} и \cite{N}.
В \S\ref{s:examp} работы \cite{AB, Be13, Ko, LV, Mo16, Ni08, Ru10, Za20} относятся к научно-исследовательской номинации, работы \cite{AB, Do, Ku, Ib, Ne, Ru10} --- к учебно-исследовательской номинации,
работы \cite{Ko, LV} --- к номинации исследовательские разработок,
работа \cite{Lo} --- к номинации наглядных материалов.
Разные версии одной работы могут относиться к разным номинациям.

\textbf{Самостоятельность} означает, что автор решил задачу сам (а не прочитал в книге и не скачал из интернета).
Конечно, сама последовательность задач может быть подсказкой; допускаются также подсказки от руководителя.
Важно, чтобы эти подсказки были полезны для развития автора, и чтобы заимствования было явно указаны в тексте.
%(например: эту лемму я доказал с существенной подсказкой руководителя, или: доказательство этой теоремы взято оттуда-то и приводится в тексте для полноты).
%Поскольку эти критерии самостоятельности не очень четкие, напишем, как мы ее проверяем.
В совместной работе должен быть указан вклад автора, с подтверждением от соавторов.
Уровень самостоятельности проясняется на Конференции, когда слушатели задают автору вопросы по содержанию работы.
%[формулировке и доказательству, о понятиях и теоремах, использованных в работе, и связях между ними владения темой]
%Если жюри сочтет, что работа не самостоятельна, то она не будет награждена.

\section{Примеры работ}\label{s:examp}

%(с пояснениями)

%\medskip
\cite{AB} Д. Ахтямов. Разрешимость уравнений третьей и четвертой степени при помощи одного радикала. 2013-2015. 
Руководители И. Богданов и А. Скопенков. 

Автору помогло знакомство с основными идеями теории Галуа (но не с ее абстрактным изложением), полученное на кружке.  Позже выяснилось, что эта задача изучалась и решена в компьютерной науке в 2002 г. (Chu, Kang).
Интересно, что Chu и Kang, ссылающиеся на абстрактное изложение теории Галуа, в одном месте применяют громоздкие вычисления вместо идей теории Галуа.
В 2013 г. работа была награждена по номинации \emph{учебно-исследовательских работ}.
После этого автор продолжил работу и в  2014 г. подготовил новую версию.
Она выложена в Архив и награждена уже научной премией.
Работа вошла в проект для ЛКТГ \cite{ABG+}.
В  2015 г. автор получил новые результаты, совместно с И. Богдановым.
Вклад автора (подтвержденный соавтором) в эти новые результаты награжден по номинации
\emph{научно-исследовательских работ}.

\medskip
\cite{Be13} В. Белоусов. A smaller counterexample to the Lando conjecture. 2013. Руководитель А. Скопенков. 

Рецензент нашел ряд ошибок в доказательстве.
Автор согласился с частью замечаний и сделал соответствующие исправления, а по другим замечаниям написал рецензенту письмо с пояснениями.
Рецензент согласился с автором и дал положительный отзыв на следующую версию работы.
Работа выложена в Архив и награждена по номинации научно-исследовательских работ.

\medskip
\cite{Do} П. Долгирев. О конкурентности некоторых чевиан треугольника. 2009. Руководитель А. Привалов. 

Первая версия работы содержала существенные ошибки и отличалась неудачным стилем изложения.
Тем не менее, совместные усилия автора, руководителя и консультанта позволили довести работу до завершения, принятия на Конференцию и награждения.
В следующем году подготовка новой работы 
%автора
потребовала значительно меньших усилий.

\medskip
\cite{Ko} С. Комаров. Две равные окружности в прямоугольном треугольнике. 2022, 2023. Автор не указал руководителя. 
%?

Автор самостоятельно развил тему, начатую публикацией Ю.А. Блинкова в <<Кванте>>.
Доведение работы потребовало более года.
Поэтому в 2022 на Конференцию была принята только короткая четкая формулировка, в номинацию \emph{исследовательских разработок}.
После этого автор продолжил работу и в следующем году подготовил полную версию.
Она выложена на геометрический сайт и награждена уже по номинации
\emph{научно-исследовательских работ}.

\medskip
\cite{Ku} А. Куликова. Теорема об изогоналях. 2016. Руководитель Д. Прокопенко. 

Задача поставлена руководителем и успешно доведена до завершения, принятия на Конференцию и награждения ней.
После Конференции работа была продолжена, ее существенно переработанная версия опубликована \cite{KP}.

\medskip
\cite{Ib} Д. Ибрагимов. 32 буквы русского алфавита кодируются последовательностями из нулей и единиц... 2015. Руководитель А. Скопенков.

Автор записал решение задачи, придуманное им на кружке.
После внесения исправлений работа была принята на Конференцию и награждена.
%по номинации учебно-исследовательских работ.
%Это характерная

%\begin{figure}[h]
%\centerline{}
%\includegraphics[scale=.7]{works2018/loshakov2.pdf}
%\caption{\cite{Lo}; пока сюда не включен, ибо 12 Мб}
%\label{f:intcubes}
%\end{figure}

\medskip
\cite{Lo} М. Лошаков.  Объединение двух пересекающихся кубов. 2018. [Рис. 
%\ref{f:intcubes} 
пока сюда не включен, ибо 12 Мб.] Руководители А. Скопенков, А. Спивак, Е. Стефанова. 

\medskip
\cite{LV} И. Васенов. Разрезания правильного многоугольника на подобные прямоугольные треугольники. 2020, 2021. Руководители М. Лацкович и А. Скопенков.

Задача поставлена руководителем по мотивам статей из научных журналов.
Автору потребовалось продолжительное время, чтобы самостоятельно реализовать идею, предложенную руководителем.
Крупный специалист в данной области, после присылки ему работы, смог получить более сильный результат.
Совместная работа автора с ним выложена в Архив и опубликована в научном журнале \cite{LV}.
Вклад автора (подтвержденный справкой от соавтора) награжден  по номинации научно-исследовательских работ.

\medskip
\cite{Mo16} E. Морозов. Обобщенная задача Аполлония. 2016-2017. Руководитель Ф. Нилов.

Автор самостоятельно исследовал обобщение на четыре окружности результата о трех окружностях, доказанного (как считается) Аполлонием более 2000 лет назад.
Работа выложена в Архив и награждена в 2016 г. по номинации научно-исследовательских работ.
После этого автор продолжил работу и получил ряд новых результатов.
Они награждены в 2017 г. по номинации научно-исследовательских работ.
Работа опубликована.

\medskip
\cite{Ne} В. Немычникова.  Об одной конике, связанной с $k$-трисами треугольника. 2014. Руководитель А. Привалов.

Первая версия работы \cite{Ne1} содержала неверную формулировку основного результата, что и было отмечено в рецензии \cite{Ne1re}.
Во второй версии \cite{Ne2} автор учла это замечание, привела правильную формулировку, а также улучшила доказательство.
Однако, формулировка следующего результата оказалась некорректной из-за отсутствия необходимых определений
\cite{Ne2re}.
Эти недочеты были исправлены в третьей версии \cite{Ne}.
Она была принята и награждена первой премией в номинации учебно-исследовательских работ.
Работа содержала элемент новизны: построенная в ней точка включена в энциклопедию центров треугольника.
Но работа не была включена в номинацию научно-исследовательских работ, потому что не было выполнено условие своевременного выкладывания в Архив.

%В первой версии работы отсутствовали четкие формулировки и завершенные доказательства.
%После обсуждения замечаний рецензента с консультантом б\'ольшая часть замечаний была учтена.
%Отзыв на вторую версию содержал лишь несколько конкретных рекомендаций по улучшению текста.
%После внесения соответствующих исправлений работа была принята и награждена по номинации учебно-исследовательских работ (поскольку новизна не была проверена выкладыванием в архив).

\medskip
\cite{Ni08} Ф. Нилов. Параболические четырехугольники. 2008. Руководитель А. Заславский. 
%М. Волчкевич.

Автор сам придумал результат.
Это получилось во многом благодаря урокам геометрии М. Волчкевича во Второй школе.
Доведение работы до принятия потребовало помощи консультанта.
Работа выложена в Архив, награждена по номинации научно-исследовательских работ и опубликована.
Впоследствии автор многократно был руководителем работ и стал постоянным членом жюри Конференции.

\medskip
\cite{Ru10} А. Рухович. На какие части разбиваются многогранники их пересечением? 2010, 2011. Руководитель А. Скопенков.

Задача была поставлена руководителем по мотивам статей из научных журналов (см. ссылки в \cite{ABR+, Ru10}).
Автор самостоятельно нашел ответ.
Работа выложена в Архив, награждена по номинации научно-исследовательских работ и опубликована в научном журнале \cite{Ru10}.
Она вошла в проект ЛКТГ \cite{ABR+}.

\medskip
\cite{Za20} Д. Захаров. О раскраске трехэлементных подмножеств. 2018. Руководитель А. Райгородский.
%?

Автор представил работу, уже выложенную в Архив и принятую в научный журнал.
Она награждена по номинации научно-исследовательских работ.

\medskip
\cite{Ze5} В. Зеленин, Идеалы в множестве квазитропических матриц. 2025. Руководитель В. Борзенко. 
 
В процессе консультаций появились новые результаты: сравните версию \cite{Ze1}, поданную на конференцию, и принятую на конференцию версию \cite{Ze5}.
Она награждена по номинации учебно-исследовательских работ. 
По-видимому, результаты работы новые. 
Поскольку эта новизна не была проверена своевременным выкладыванием в Архив, работа не была включена в номинацию научно-исследовательских работ, 

\medskip
\cite{Zu6} К. Зюбин. Минимальное время получения натурального числа $n$ из 1 путем выполнения заданных операций. 2020. 
%??? Руководитель А. Привалов. 

В процессе консультаций были упрощены и прояснены формулировки: сравните версию \cite{Zu1}, поданную на конференцию, и последующую версию \cite{Zu3}.   
При этом на конференцию была принята версия \cite{Zu6}, награжденная по номинации учебно-исследовательских работ.

\medskip
Тематика работ разнообразна, больше всего работ подается по элементарной геометрии и дискретной математике.  %(комбинаторике?).

В вышеприведенных примерах \cite{Ni08, Ru10, Mo16, AB, LV} авторам потребовалось 1-2 года, чтобы довести доказательство до завершенного.

%?Стандартные ситуации (отклонение, задачи - уч-иссл, мутные тексты - иссл. разработки)

\section{Завершенность (надежность) доказательств}\label{zaver}
 
\textbf{Зачем?}

Важнейший результат интеллектуальной деятельности --- текст (или программа), который можно прочитать (или которой можно воспользоваться) без консультаций с ее автором.
И который проверен его авторами настолько серьезно, что пользование не принесет вреда.
Поэтому \emph{важное отличие исследовательской деятельности от  учебной --- необходимость серьезно проверять результаты.}

При обучении развивается \emph{способность} что-то сделать, обычно связанная с \emph{имитацией} деятельности для пользователя 
(построить модель самолета, решить задачу на олимпиаде).
Продвинутое обучение предполагает постепенное развитие в сторону \emph{реализации этой способности результатом} 
%[этой способности в результативную] 
(починить водопроводный кран в настоящем самолете, опубликовать короткое завершенное доказательство олимпиадной / кружковской задачи, даже не претендующее на новизну).
Ученик, прошедший такое обучение,
%даже если результат теперешней работы ученика известен в науке,
в будущем сможет писать полезные и надежные научные работы и программы, проектировать атомные станции и разрабатывать новые лекарства.
%другие экономические, индустриальные и медицинские объекты.
См. также \cite[абзац на стр. 24-25]{Sk08} (о парадоксе <<отрицательного вклада>>).

%Есть доказательства, <<правильные>> только благодаря наличию учителя-оракула, т.е. человека, знакомого с фактом, его доказательствами и тем, что вокруг.
%Если же сам автор --- без проверки оракулом --- будет писать доказательства такого уровня, то они часто будут ошибочными.

%(Конечно, не обязательно серьезно проверять результаты, пользование которыми не предполагается.
%[Часто те, кто ленятся серьезно проверять результаты, задумываются не об их использовании, а о возможности получить награду без проверки.])

%Поэтому четкие формулировки и завершенные доказательства считаются необходимыми (кроме новизны) для публикации математического результата в рецензируемом научном журнале.
%Четкие формулировки результатов необходимы для их использования.

В научном мире подготовка надежных текстов для пользователя происходит в рамках системы <<тщательного рецензирования>> (peer review).
%https://en.wikipedia.org/wiki/Peer_review
%https://en.wikipedia.org/wiki/Peer_critique
%https://en.wikipedia.org/wiki/Scholarly_peer_review
При ней окончательный (публикуемый, награждаемый) текст формируется в процессе работы автора над замечаниями рецензента.
Этот процесс аналогичен процессу \emph{отладки} программы, являющемуся основным в программировании.
Учителю и (постепенно) ученику полезно знакомиться с таким стилем работы и приобщаться к нему.

%Отсутствие четкой связи
%(обратите внимание на парадокс отрицательного вклада вверху стр. 25 и в замечании 1.5, соответственно).

%?
Напомним также (см. введение), что написание завершенных доказательств способствует развитию интуиции за пределы данного материала.

\textbf{Что такое завершенное (надежное) доказательство?}

В нашем понимании это --- доказательство, вероятность ошибки в котором ничтожно мала благодаря тому, что оно удовлетворяет высоким стандартам
строгости и структурированности.
Для надежности (и тесно связанной с ней  понятности) математического текста человечество выработало некоторые стандарты проверяемости.
Только текст, в котором эти стандарты соблюдаются, мы называем <<завершенным доказательством>>; соответствующий термин в науке --- <<надежный источник>> (<<reliable reference>>).
Другими словами, оно написано уже для пользователя, а не еще для разработчика (work in progress) и не еще для проверяющего (как на олимпиаде).
(Завершенность понимается именно в этом смысле, а не в смысле отсутствия направлений для дальнейших исследований.)

Завершенность доказательства в научно-исследовательской номинации
%а в также наиболее грамотных и содержательных работах учебно-исследовательской номинации, нет, ибо научный
соответствует принятию работы в ее текущем виде в рецензируемый научный журнал.
% (с доработками на усмотрение автора).
Завершенность доказательства в наиболее грамотных и содержательных работах учебно-исследовательской номинации понимается аналогично.

Таким образом, школьник, учась писать завершенные доказательства на простых коротких примерах, серьезно облегчает себе написание завершенных доказательств в последующих научных работах (или выполнение работы для пользователя вне науки).

\textbf{Примеры.} 
%[См. сноску \ref{f:examples}]

Приведем яркую цитату из \cite{JQ93}:\footnote{Многие ситуации, связанные с нестрогостью в математике, были менее удачными, чем описанные выше. 
В особенности это проявляется там, где нестрогий материал подается как известный и надежный, с претензией на приоритет. 
Иногда это <<честная ошибка>>, иногда результат нестандартного понимания того, что является доказательством. 
\newline
Прямые ошибки приносят меньше вреда. 
Например, фундаментальная лемма Дена о двумерных дисках в трехмерном многообразии была представлена в 1910 г. 
В доказательстве была обнаружена ошибка, и ко времени появления верного доказательства (К. Папакириякопулос, 1957) эта лемма была признана важной гипотезой. 
\newline
Низкие стандарты доказательств приводят к более значительным трудностям. 
В восемнадцатом веке небрежные рассуждения привели к болезненным проблемам в области сходимости последовательностей и равномерной сходимости функций. 
В качестве <<лечения>> была введена строгость. 
Она была вовремя принята, несмотря на возражения ряда математиков. 
Это помогло избежать более значительных проблем.   
\newline
Позднее в двадцатом веке значительных проблем не избежала итальянская школа алгебраической геометрии. 
Она пришла в упадок после поколения блестящих нестрогих рассуждений. 
Обсуждение трудностей и долгого оздоровления приводится в [EH, K].
В 1946 к этому предмету по-прежнему относились с таким подозрением, что А. Вейлю пришлось защищать свой интерес к нему; см. введение к [W].}

\emph{Most of the experiences with theoretical mathematics have been less positive
than those described above. This has been particularly true when incorrect or
speculative material is presented as known and reliable, and credit is claimed by
the perpetrator. Sometimes this is an `honest mistake', sometimes the result of
nonstandard conceptions of what constitutes proof. Straightforward mistakes are
less harmful. For example, the fundamental `Dehn lemma' on two-disks in three-manifolds was presented in 1910. An error was found\footnote{П.С. Александровым, который считал нахождение этой ошибки своим важным достижением (прим. авторов).}, and by the time it was proved (by Papakyriakopolos in 1957) it was recognized as an important conjecture.}

\emph{Weak standards of proof cause more difficulty.
In the eighteenth century, casual reasoning led to a plague of problems in analysis concerning issues like convergence of series and uniform convergence of functions.
Rigor was introduced as the antidote.
It was adopted over the objections of some theorists in time to avoid major damage.}

\emph{More recently in this century the `Italian school' of algebraic geometry did not
avoid major damage: it collapsed after a generation of brilliant speculation.
See [EH, K] for discussions of the difficulties and the long recovery.
In 1946 the subject was still regarded with such suspicion that Weil felt he had to defend his interest
in it; see the introduction to [W].} 
%\cite{JQ93}

Применение сформулированных принципов рецензирования (а также умение сделать соответствующую работу автора интересной в процессе консультаций) хорошо прослеживается на примере работы В. Немычниковой в \S\ref{s:examp}. 
%???(руководитель А.А. Привалов), см. .
%? Развитие профессионализма в рамках участия в ММКШ (а также, и даже в первую очередь, в рамках обучения в матклассе) помогло В. Немычниковой успешно учиться в выбранном ВУЗе (МФТИ).

Более подробные рассуждения и конкретные примеры см. в \cite{JQ93, ABC+, BBB+, Sk08, Sk21d}.

\newpage
\textbf{Как конкретно ученик подходит к завершенности доказательств?}

Обучаясь в <<системе кружков и олимпиад>>, школьники постепенно переходят от занимательных задач к математическим   теориям (ср. \S\ref{why}).
%ГЧ в большинстве случаев нет необходимости в тщательной проверке решений, прием задач отчасти превращается в "читатель ждет уж рифмы "розы" -- так на, возьми ее скорей". Отсутствие небходимости делает несколько менее удобной тренировку этого навыка. В ситуации же исследовательских задач
%?Правильно не избегать естественного нарастания
Нарастают необходимость и возможность контролировать правильность (и понятность) решений --- особенно, когда обучение включает элемент исследования с неизвестными заранее ответом или подводными камнями в доказательстве.\footnote{К сожалению, в погоне за количеством в ущерб качеству такая необходимость не всегда признается, а возможность --- реализуется.}
%, особенно в ВУЗах.)
%[есть распространенное заблуждение] 
В частности, необходимость записи \emph{завершенных доказательств} некоторых утверждений.
А значит, и \emph{четких формулировок} утверждений (например, промежуточных шагов решения).

К завершенности некоторых доказательств стоит стремиться только на достаточно продвинутом этапе.
Тем более важно не пренебрегать таким стремлением, особенно когда ученика за него можно наградить.
Интересно, что обычно успешные олимпиадники достаточно мотивированы и подготовлены к переходу на этот более высокий уровень.

%Кроме того, приведение четких формулировок результатов и промежуточных шагов решения (например, лемм) является частью написания завершенных доказательств.

После \emph{неформальной проверки} рекомендуем проверять доказательство \emph{записью}. 
Разумно также пользоваться помощью 
%учителя, руководителя, 
консультанта 
%математика близкой специальности, 
и даже однокашника.

Разумеется, научить писать завершенные доказательства сможет только тот, кто сам это умеет и от себя требует.

%[ГЧ Вот этот фрагмент я бы очень хотел видеть в финальном варианте текста (причем желательно в жестком варианте), но для этого к нему должна быть сделана некоторая подготовительная работа.]
%Учителя / руководители могут научить писать завершенные доказательства, если сами неплохо умеют делать это и учатся делать лучше (либо прошли такую <<школу>> в прошлом).

В рамках кружков начальное знакомство с системой <<тщательного рецензирования>> фактически реализуется \emph{письменной} сдачей задач (такое было и на ранних ЛКТГ).

Мы знакомим публику с этой системой,
% <<тщательного рецензирования>>,
выкладывая \cite{R} все работы, поданные на Конференцию, и рецензии на них.
Ведь конкретные примеры рецензирования гораздо более эффективны, чем общие слова (см., в частности,
%характерные
примеры \cite{Be13, Ne} в \S\ref{s:examp}).
% или обращайтесь за консультациями!

%Получаемые версии текста делают возможными все более и более серьезные проверки доказательства --- как ими, так и более далекими специалистами, уже не знакомыми с устным изложением работы.

%Писать четкие формулировок и завершенные доказательства разумно до направления работы в научный журнал или на награду, предполагающую их наличие.
%Только тексты, в которых эта работа уже проделана, прилично присылать на проверку квалифицированным рецензентам (обычно людям занятым).

%Например,
Опишем
%более продвинутую
реализацию такого знакомства
%с системой <<тщательного рецензирования>>
при представлении некоторых тем на ЛКТГ 2010х-20х годов.
Команда, решающая данный цикл задач,
%задачи по данному тексту,
может состоять из любого числа участников.
За каждое решение, \textit{написанное для пользователя}, оцененное в <<+>> или <<+.>>, команда получает пять <<бобов>>  (см. рекомендации \cite{Re}).  %\url{https://www.mccme.ru/circles/oim/home/pism.pdf}).
Дополнительные бобы могут выдаваться за красивые решения, решения сложных задач или оформление некоторых решений в системе \TeX.
У жюри бесконечно много бобов.
У каждой команды в начале десять бобов.
Если у команды есть бобы, то решения можно сдавать и \textit{устно}, и как \textit{написанные для соавтора} (т.е. менее точно, чем для пользователя).
Команда отдает один боб за каждую такую попытку сдачи (удачную или нет).
Участники сами выбирают задачу, по которой пишут решение для пользователя, но могут обсудить с жюри свой выбор.
Приглашаем участников обсуждать с жюри возникающие вопросы (без траты бобов).
В частности, жюри может помочь записать решение для пользователя или дать подсказку к задаче, которая долго не решается.
%На Летних Конференциях выдается часть данного текста.
Те, кто успешно сдают задачи, завоевывают право получить интересные \emph{дополнительные задачи для исследования}.
% (из не выданной части).

\section{Чем мотивировано различение по номинациям?}\label{s:nomi}

\textbf{Работы, не претендующие на научную новизну.}
%\footnote{[В (неформальных) обсуждениях ответственный секретарь редколлегии Мат. Просвещения (МП) предложил удалить данный кусочек, \S\ref{s:howto} и б\'ольшую часть \S\ref{s:rest}. 
%Предложения о сокращении как таковые мы приветствуем. 
%Более того, 
%\newline
%$\bullet$ перед представлением статьи в МП мы сделали огромную работу по сокращениям; в результате обсуждений с ответственным секретарем двух из этих сокращений в одном случае он предложил вернуть сокращенный текст; 
%\newline
%$\bullet$ АС сразу при представлении других статей в МП указывал на возможность удалений отдельных параграфов; никогда такие возможности редколлегией МП не %одобрялись.  
%\newline
%Но, по нашему мнению, удаление указанных кусочков покажется неуместным большинству читателей МП (<<между Квантом и научными журналами>> по словам ответственного секретаря). 
%Это мнение основано на обширных профессиональных обсуждениях с учителями, студентами и школьниками на протяжении десятилетий. 
%Однако мы не претендуем на истину в последней инстанции, а предлагаем явно спросить мнение учителей, 
%студентов и школьников о желательности предложенного удаления кусочков. 
%При этом конкретные сокращения в этих кусочках мы готовы делать (или дать кусочки мелким шрифтом, etc).]}

Целью исследовательской деятельности школьника не обязательно должна быть завершенная научная работа.
Большинству школьников трудно ее сделать.
Тем более, что начинающему не очевидно, что имеются разные степени завершенности: от work in progress (на стадии написания первых версий или текста доклада на научном семинаре; похвала руководителя и коллег) до reliable reference (опубликованная работа, которой может пользоваться другой математик или инженер; научная премия или возможность защитить курсовую, дипломную, диссертацию и т.д.).

Cамостоятельное решение важной красивой (известной) задачи или проведение вычислительного эксперимента как раз подведет школьника и к разумному выбору задач для более серьезной работы, и к необходимости четких формулировок и завершенных доказательств (\S\ref{zaver}), и к возможности их написать.
Награждены могут быть все работы --- при условии, что известность результата (или отсутствие серьезной
проверки новизны) явно отражена, а итоги экспериментов или незавершенных доказательств не называются
теоремами.
Введенное для этого явное различение по номинациям поможет школьникам

$\bullet$ решать доступные и полезные им задачи;

$\bullet$ учиться серьезно проверять свои доказательства;

$\bullet$ сформировать реалистичное представление о научной работе.

\textbf{Научно-исследовательские работы.}

Мы считаем, что работы школьников, награждаемые премиями по этой номинации, должны удовлетворять критериям
(полноты доказательств и серьезности проверки новизны), предъявляемым к \emph{завершенным научным работам}
(без скидки на возраст авторов).
При этом требования к глубине и количеству результатов могут быть более слабыми.

По нашему мнению, сейчас некорректно утверждать новизну результатов работы без ее подачи в архив \url{http://arxiv.org} или \url{https://hal.science} и отсутствия 
в течение нескольких недель информации об известности этих результатов.
%[(в течение нескольких недель после появления) писем  автору со ссылками на те результаты работы, которые уже известны].
Такая проверка новизны имеет свои ограничения
%(реакция сообщества возможна, но маловероятна)
и не отменяет важнейшей предварительной работы (рассылку специалистам по данной тематике и поиск по математическим базам данных).
Однако она проста для включения в формальные требования, а также совмещает проверку новизны с пользой для сообщества (работа доступна) и 
ответственностью рекомендателя / автора (удалить уже выложенный текст, как правило, невозможно).

Выложить первую работу на архив невозможно без рекомендации (математика, имеющего такое право).  
%достаточное количество работ по данной теме в архиве за последние годы).
При попытке выложить работу в Архив автору предлагаются <<предостережения>> (ссылка на которые продублирована на сайте Конференции).
Рекомендации помогает получить программный комитет,  
%и приглашенные им потенциальные рекомендатели.
 если (по его мнению) выкладывание работы не нанесет ущерб репутации ее автора. 
% выкладывание имеющейся версии работы не создаст отрицательного мнения об уровне автора
Получению рекомендации предшествуют содержательные обсуждения: перед выкладыванием часто автор исправляет недочеты, иногда выясняется известность результата и / или возможность его развития (см. пример работы К. Зеленина в \S\ref{s:examp}).
Обычно в результате мы даем рекомендацию, а также помогаем автору выложить представленную в Архив версию
%со штампом архива submit*.pdf
на популярный ресурс по данной теме (например, на \url{https://geometry.ru}) --- на случай задержки или отклонения модераторами Архива.
% (но если к тексту имеются серьезные замечания, то предлагаем автору обратиться к ).
Некоторые работы выложены (и могут быть выложены) в Архив независимо от программного комитета --- на их рекомендателях и лежит ответственность за рекомендацию.

\textbf{Пример.}\footnote{\label{f:examples} Этот пример и примеры из \S\ref{zaver} добавлены в соответствии с замечанием редколлегии Мат. Просвещения о недостаточном количестве конкретных примеров.} 
%(а не только ссылок на примеры; такие ссылки уже были).  
%Мы согласны с ответственным секретарем, что одни примеры можно (обратно) заменить на ссылки, а другие удалить. 
%Мы будем рады это сделать, если редколлегия подтвердит, что соответствующее замечание снято (даже при отсутствии явного приведения примеров вместо ссылок) и не будет далее повторяться.] 
%Следующий пример не столь значителен для математики, зато уже относится к конференциям школьников.

Работа С. Тищенко (руководитель А.Я. Канель-Белов) была награждена 1й премией на международной конференции школьников Intel-ISEF в 2000 г.
Работа претендовала на решение одной из проблем Эрдеша.
За несколько месяцев до этого руководителю была направлена (Б. Мохаром через А. Скопенкова) ссылка на решение этой проблемы, полученное ранее.
Однако о том, что эта проблема решена ранее, не сообщалось ни в работе, представленной на Intel-ISEF, ни  в ее опубликованной версии \cite{Ti01} (в аннотации и истории вопроса на первых страницах).\footnote{В конце статьи \cite{Ti01} есть теорема 1а, названная отрицательным решением проблемы Эрдеша, а еще через страницу дана ссылка на работу, в которой была получена теорема 1a.
Эта информация должна присутствовать при первом упоминании в статье решения проблемы Эрдеша.}
%?(возможно, вместе с информацией, что )
Следующая информация, возможно, не была тогда сообщена руководителю и автору, но неизбежно возникла бы в работе при добросовестном написании истории вопроса.
The paper
%given under the provided link %[http://www.mathnet.ru/php/archive.phtml?wshow=paper&jrnid=fpm&paperid=557&option_lang=eng]
solves a very specific problem of diameter 3 and subcubic graph, which I doubt can be referred to as a problem of Erdos (the general degree-diameter problem is related to Erdos, though).
In particular, there is no reference to a paper by Erdos.
%Mohar

%в случае принятия в МП мы обязуемся в течение 10 дней перевести все английские кусочки на русский язык.

На ММКШ такое практически невозможно ввиду прозрачности принятия в номинацию научной-исследовательских работ премией (в частности, выкладывания в архив).

\section{Рецензирование и  консультации (заочный первый тур Конференции)}\label{s:rest}

В этом и следующем параграфах мы переходим к наиболее существенной --- практической --- и нетривиальной части статьи. 
%!!!
Трудность не в высказывании позитивных принципов, а в следовании им. 
Характерная реакция многих математиков и учителей на начало данной статьи --- <<все это замечательно, но невозможно реализовать>>. 
Поэтому приведенное описание реализации представляется нам важным.

%\small
{\bf Как устроено рецензирование?}

Поданные работы выкладываются на страницу рецензирования \cite{R} (зачем? см. ниже о прозрачности).
%Они направляются на \emph{рецензию}.
Работы рецензируются для \emph{проверки} результатов и их доказательств.
Первый шаг такой проверки --- проверка формулировки (является ли она четкой? нет ли очевидных контрпримеров?).
%В олимпиадах участнику не приходится проверять четкость формулировок задач, 
Однако умение четко сформулировать лемму (в частности, попытаться построить контрпример к придуманной формулировке) помогает в решении сложной олимпиадной задачи.

Взглянуть на свою работу глазами другого человека, <<пользователя>>, трудно.
Однако именно взгляд со стороны (рецензента и консультанта) помогает увидеть и ошибки, и новые возможности, которые трудно заметить автору или руководителю, <<живущим внутри>> задачи.
Без обсуждений с математиком, не знакомым заранее с работой, автору обычно невозможно подготовить по-настоящему проверенный доклад.
%, интересный слушателям (тоже изначально не знакомым с работой).

%Многие работы, принимаемые на Конференцию, не претендуют на новизну (результатов или доказательств) и не содержательны вне пределов <<кружковской>> задачи (\S\ref{s:kaki}).
%%Поэтому рецензент не пишет о новизне или высокой содержательности. Однако
%По работам, имеющим новизну или высокую содержательность, рецензия оценивает степень новизны или содержательности.
Рецензент лишь оценивает представленный текст, но не должен учить автора писать тексты (в отличие от руководителя и консультанта).

{\bf Ответственность 
%рецензента и 
программного комитета.}

Рекомендации рецензента (принять / отклонить работу) обосновываются (в самой рецензии, а не в телефонных разговорах с программным комитетом) 
конкретными замечаниями со ссылками на конкретные места в работе.
Рецензия выкладывается на сайт рецензирования без указания фамилии рецензента.
Кроме того, любой желающий может прислать \emph{подписанную} рецензию на любую работу.
Такие рецензии также выкладываются на страницу рецензирования (даже те, с которыми программный комитет не
согласен, кроме спама).
Поэтому программный комитет

%$\bullet$\footnote{[По желанию редколлегии МП, если первые две жирные точки не ясны, то они могут быть удалены.]}
% рецензент, приглашенный программным комитетом, несет ответственность за обоснованность рецензии перед ним (ответственность рецензента обеспечивается репутационным риском вследствие некомпетентного отзыва);

%$\bullet$ рецензент, приславший рецензию, несет ответственность за ее обоснованность рецензии перед сообществом;

$\bullet$ несет ответственность перед сообществом за обоснованность решения о принятии / отклонении работы на основании рецензий ;

$\bullet$ соблюдает принцип непредвзятости: оценивается именно работа, а не личность автора;

%decisions ought to be based on objective criteria, rather than on bias, prejudice, or preferring to benefit one person over another for improper reasons;

$\bullet$ сооблюдает принцип <<решение, принимаемое закрытым обсуждением, не должно  противоречить публичным обсуждениям содержания работы>>.
\footnote{Например, работа не может быть награждена премией \emph{второй} степени (по какой-то номинации) на закрытом заседании жюри, если опубликован единственный отзыв, рекомендующий работу к \emph{высшей} награде (по этой номинации), и на публичном обсуждении доклада не было высказываний по содержанию работы, ограничивающих сверху ее оценку  (например, <<имеются такие-то недочеты>> или <<новизны недостаточно для публикации>>).
Благодаря этому принципу автор (и все участники публичного обсуждения) могут ответить на все высказывания по \emph{содержанию} работы, влияющие на \emph{закрытое} обсуждение ее \emph{оценки}.
Заметим также, что без выполнения такого принципа невозможно утверждать, что выполняется принцип непредвзятости.}

%Confidential decisions should match public discussions.
%Every substantial scientific argument affecting a confidential decision should be publicly available.
%Here a `decision' is a decision by an editorial board, dissertaion/diploma committee, etc.
%An example on dissertation defence in Russia: if there are no public objections to a dissertation, but a dissertation committee votes against the dissertation, then the committee должен быть is dissolved.
%The {\it quantity} of persons publicly supporting some point of view need not be matched by a confidential decision.

О такой системе <<прозрачного тщательного рецензирования>> (transparent peer review) см. подробнее \cite{T}.

{\bf Своевременность.} 

Рецензия готовится в течение не более чем трех недель после получения работы.

{\bf Что делать, если прорецензированная работа отклонена?} 

Автор может подать новую версию, учитывающую замечания рецензента.
Либо новая версия принимается, либо процесс повторяется (возможно, несколько раз).
Этот процесс можно ускорить, упростить и сделать более приятным, используя перед подачей очередной версии
рекомендации (\S\ref{s:howto}) и  консультации (см. ниже).
Иногда этот процесс не завершается до крайнего срока в данном году, тогда он переносится на следующий год.
Он аналогичен работе научных журналов.
Он отличается от олимпиад (и контрольных) тем, что время у автора на доработку не ограничено, но принимается только доработанная (достаточно надежная) версия, а никаких <<баллов>> за недоработанную версию не ставится.

\emph{Несогласие автора с замечанием рецензента} --- нормальная, хотя и редкая ситуация \cite{R}.
Вот способы ее разрешения.

(1) Автор может на консультации обсудить те замечания рецензента, которых не понимает или с которыми не согласен.
Если нужно, консультант поможет автору грамотно составить письмо рецензенту, см. (2).
Эти обсуждения \emph{не выкладываются} на страницу рецензирования.

(2) Автор может прислать письмо о замечаниях рецензента
%, с которыми он не согласен,
с обоснованием несогласия.
Это письмо \emph{выкладывается} на страницу рецензирования.
Поэтому автору желательно обсудить текст письма с руководителем и с консультантом: ведь текст письма, как и текст работы, показывает уровень математической грамотности автора.
%Рецензент ответит на письмо не позже, чем через три недели.

{\bf Прозрачность.}

Выкладывание работ в интернет позволяет избежать {\it профанации}, т.е. несоответствия между заявленным (скажем, на сайте или в письме к спонсору) уровнем конференции-конкурса и реальным уровнем работ (который может даже признаваться в частных разговорах).
Авторы учатся нести репутационную ответственность за качество своих текстов.
При этом
%по нашему мнению
наличие \emph{исправленной} автором ошибки обычно повышает
%(а не понижает)
репутацию автора.
%если недочеты в работе исправляются, то выкладывание на страницу резензирования версий, в которых недочеты исправляются, обычно повышает репутацию автора.

Выкладывание работ в интернет позволяет также
%[полезно и интересно]
педагогическому сообществу получить представление о системах <<тщательного рецензирования>> (\S\ref{zaver}) и <<прозрачного тщательного рецензирования>> \cite{T}.
Для такого рецензирования необходимо рассматривать работы заранее.

%\section{Консультации}\label{s:cons}

\textbf{Мы предлагаем авторам консультации.}

Мы приглашаем и руководителей участвовать в них.
Работы, присланные на консультацию, не выкладываются на страницу рецензирования.
Консультант дает замечания и советы,
%по работе над замечаниями,
а также отвечает на вопросы.
Это помогает автору подготовить версию работы, которая получит положительный отзыв рецензента и будет принята.
При этом одобрение консультанта не означает принятия работы.
% (обычно рецензент и консультант --- разные люди).

%\normalsize

%Обычно требуется несколько консультаций.
%Поэтому консультации проходят круглогодично.
%Мы приглашаем школьников присылать работы на консультацию задолго до крайнего срока подачи работ.
%Некоторые школьники, награжденные научной премией, консультировались по своим работам год и более.

\section{Рекомендации: как начинающему написать работу по математике}\label{s:howto}

%Сначала мы приведем сами рекомендации, а потом их мотивируем.

%\textbf{Что д\cyryatлать?}

%\small

Эти рекомендации не касаются важных предварительных обсуждений с учителем (руководителем),
в которых вырисовывается решение задачи (или программа или рисунок).
Мы рекомендуем, что делать \emph{после} таких обсуждений.

\textit{Если ученику пока не интересно учиться приводить четкие формулировки и завершенные доказательства} (\S\ref{zaver}), то рекомендуем подать на Конференцию работу в номинацию  наглядных / экспериментальных материалов (в которой такие формулировки и доказательства не требуются, см. \S\ref{s:kaki}).

\textit{Работа по математике} --- это в первую очередь четко сформулированное математическое
утверждение (теорема или гипотеза).
Только после четкой формулировки утверждения имеет смысл писать его доказательство и мотивировки.
Хотя \emph{придумывание} может происходить в другом порядке, даже в обратном.
См. <<Зачем>> в \S\ref{zaver}.

%и даже необходимо

\textit{Если ученик уже дошел до того уровня, когда уже интересно и полезно учиться приводить четкие формулировки и завершенные доказательства,} то полезно написать четкую формулировку \emph{одного} результата, включая строгие определения используемых в ней понятий (почему только одного? см. следующий абзац).
Если в голове или в тексте имеется несколько результатов, то выберите любой --- например, самый интересный или самый простой для записи.
Начните с удаления всего текста, кроме этой формулировки (и используемых в ней определений).
В частности, не нужно употреблять и определять понятия, не используемые в формулировке результата.
Учитывайте рекомендации \cite{Re}.
Обсуждайте формулировку с руководителем и консультантом; исправляйте ее в соответствии с их замечаниями.
Такие обсуждения гораздо быстрее, проще и приятнее, чем с рецензентом.
Когда руководитель и консультант признают формулировку четкой, подайте ее на Конференцию, назвав результат гипотезой (поскольку в тексте не будет ее доказательства).
Если рецензент признает формулировку четкой (и в работе не будет ничего другого), то работа будет принята в
номинацию исследовательских разработок.
%Можно подать на Конференцию сразу несколько формулировок; тогда для принятия работы необходима четкость <i>каждой</i> из них.
Как правило, только после такого принятия четкой формулировки полезно (с помощью руководителя и консультанта) добавлять новые куски (например, эвристические соображения или завершенные доказательства) --- если ученик сочтет это нужным.
Некоторые работы, начатые с короткой четкой формулировки, уже доведены до полных текстов, награжденных научной премией Конференции (см., например, \S\ref{s:examp}, \cite{Ko}).

\textit{Главное --- не количество, а качество!}
Минимизировать текст нужно, ибо недостаточная математическая грамотность какой-то его части может помешать принятию всей работы.
Длина текста сама по себе не является достоинством (она является достоинством только тогда, когда
текст длинный за счет глубины и нетривиальности результатов, а не за счет неудачного изложения).
\textit{Один} четко сформулированный результат, имеющий завершенное доказательство, как правило, является
заметным достижением.
Его подготовка принесет пользу автору, а его публикация --- пользователю.
Обычно такая подготовка требует нескольких итераций <<консультация (рецензия) --- новая версия>>.
Когда текст будет одобрен рецензентом, можно будет добавить введение (объяснение, почему результат интересен и как он связан с другими) и дальнейшие результаты.
Большое количество результатов, не проверенных по-настоящему, не нужно пользователю.
Поэтому работа сразу над многими результатами до того, как рецензент одобрил \textit{один} из них, обычно является пустой тратой времени и осложняет обучение.
Если некоторые рассуждения не являются завершенными доказательствами, то либо удалите их, либо назовите
их набросками, а их выводы --- гипотезами, чтобы не дезориентировать читателя (и не помешать принятию
работы на Конференцию).

См. \cite{S} и истории рецензирования работ на \cite{R} (например, К. Зюбина в 2020 и О. Кашурина в 2019).

%\normalsize

%\textit{Вопрос.} У меня катастрофически не хватает времени, потому что [перечисление].
%Я постараюсь учесть все полученные рекомендации потом.

%\textit{Ответ.} Понимаю Вас, нам всем приходится наступать на горло одним своим песням, чтобы спеть другие.
%Желательно, чтобы Вы либо подготовили до какого-то конкретного срока завершенное доказательство, либо предложили сделать это другим желающим школьникам.

%\section{Какие проблемы имеются Чем мотивированы правила Конференции?}\label{s:motiv}
%Правила конференции отражают реальную жизнь, поэтому они не так просты.
%Приведем объяснения, показывающие естественность этих правил.

\section{Руководство Конференции, его открытость}\label{s:manag}

{\bf Программный комитет:}

$\bullet$ Д.В. Трещёв (председатель), академик РАН, директор Математического Института им. В.А. Стеклова;

$\bullet$ А.А. Заславский, кандидат техн. наук, председатель жюри олимпиады им И.Ф. Шарыгина по геометрии;
%главный редактор журнала <<Journal of Classical Geometry>>;

$\bullet$ А.Б. Скопенков, доктор физ.-мат. наук, профессор МФТИ и НМУ;

$\bullet$ И.В. Ященко, кандидат физ.-мат. наук, директор МЦНМО.

{\bf Жюри:} члены программного комитета, руководители принятых работ школьников (по одному от каждой работы) и следующие постоянные члены жюри:

$\bullet$ П.В. Бибиков, кандидат физ.-мат. наук, зав. кафедрой математики лицея <<Вторая школа>>;

$\bullet$ В.В. Буланкина, учитель математики школ 2007 и 444, методист ЦПМ.

$\bullet$ Ф.К. Нилов,  учитель математики школы 57;

$\bullet$ А.А. Привалов,  учитель математики химического лицея N 1303;

$\bullet$ Д.В. Прокопенко, учитель математики школы 2007;

$\bullet$ А.М. Райгородский, доктор физ.-мат. наук, директор Физтехшколы прикладной математики и информатики МФТИ;

$\bullet$ А.И. Сгибнев, кандидат физ.-мат. наук, зав. кафедрой математики школы <<Интеллектуал>>.

{\bf Оргкомитет:} И.В. Ященко и Б.Р. Френкин.
%, кандидат физ.-мат наук, МЦНМО.

{\bf Совет рецензентов:} \cite{Rf}.

В создании конференции важную роль сыграли Н.Н. Константинов (1932-2021) и В.Д. Арнольд (1968-2017).
В программном комитете состояли также Н.Н. Константинов в 2008-2021 гг и А.И. Сгибнев в 2012-2021 гг.

{\bf Открытость.}
Руководители, работы учеников которых регулярно принимаются на Конференцию, по решению программного комитета и постоянных членов жюри приглашаются в постоянные члены жюри.
Авторы многочисленных удачных рецензий по решению программного комитета приглашаются в совет рецензентов.
Постоянные члены жюри, активно сотрудничающие с программным комитетом, по решению программного комитета приглашаются в программный комитет.

\section{Приложение: характерные заблуждения}\label{s:wrong}
 
Вред от многих заблуждений очевиден. 
Поэтому, как правило, они не формулируются явно, но используются при принятии решений. 
Мы попробовали вычленить и явно сформулировать эти заблуждения.\footnote{В этом нам помогла статья \cite{Be14}. 
Мы благодарны ее автору, Алексею Канелю-Белову, за письменную формулировку своей точки зрения, из которой легко вычленить заблуждения, см. ниже. 
Скажем очевидное --- в статье \cite{Be14} много мыслей и бесспорных, и спорных, но которые не стоит считать заблуждениями.}  

\begin{blud}\label{b:impo}
Школьник не может сделать научную работу. 
\end{blud}

Это опровергается примерами, см. в частности,  \S\ref{s:examp}. 
%\cite[\S3]{V2}. 
Иногда для опровержения заблуждения \ref{b:impo} приводятся примеры незавершенных работ школьников (уровня `work in progress', а не `reliable reference').  
Такие примеры лишь поддерживают заблуждение \ref{b:impo}. 

\begin{blud}\label{b:easy}
Школьнику легко сделать научную работу. Многие школьники делают научные работы.  
\end{blud}

Подготовке научной работы предшествуют годы изучения <<внешкольной>> математики. 
Характерный ориентировочный срок подготовки первой короткой научной работы --- полгода-год. 
Наиболее удачна ситуация, когда учение настолько постепенно и естественно перерастает в исследование, что одно трудно (и не нужно) отделять от другого. 
Тогда невозможно определить, где кончается учебная работа и начинается научная. 
Обсуждаемый срок зависит от доли времени, которое начинающий посвящает исследованиям --- оно и не должно быть большим по сравнению со временем на обучение. 

Одна из трудностей в подготовке первой научной работы начинающим --- осознание разницы между уровнем `work in progress' и `reliable reference', 
и переход на более высокий уровень. 
В (не)преодолении этой трудности велика роль руководителя (и окружающих вообще). 
Распространение понимания такой разница и путей перехода на более высокий уровень способствует увеличению количества научных работ начинающих (и их гармоничному развитию) 
На \cite{R} можно посмотреть, как реально `work in progress' превращается в завершенную научную работу, или награждается как таковой 
(без обмана автора). 

\begin{blud}\label{b:arxiv}
Организаторы конференций, ставящие --- в качестве необходимого yсловия признания награждаемой работы научной --- ее выкладывание в архив, наносят вред авторам работ. 

Вот формулировка этого заблуждения, которая 
%менее ясна, зато 
опубликована уважаемым автором в уважаемом журнале \cite[(1) в п. 3.1 на стр. 240]{Be14}.     

<<\emph{Не надо учащихся <<подставлять>>. Например, сам учащийся может неадекватно оценивать свою
работу и её изложение, так что требование выкладывать работу на сайте www.arxiv.org в качестве предварительного условия окажется провокационным.}>>
\end{blud}

Это заблуждение основано на заведомо ложном предположении о том, что любой желающий может выложить в Архив свою первую работу.  
См. по этому поводу абзац перед кусочком <<пример>> в \S\ref{s:nomi}.  
%в \cite[абзац перед кусочком <<пример>> в \S5]{V2}: 
%`\emph{Выложить первую работу на архив невозможно без рекомендации (математика, имеющего такое право).  
%При попытке выложить работу в Архив автору предлагаются <<предостережения>> (ссылка на которые продублирована на сайте Конференции).
%Рекомендации помогает получить программный комитет, если (по его мнению) выкладывание работы не нанесет ущерб репутации ее автора. 
%Получению рекомендации предшествуют содержательные обсуждения: перед выкладыванием часто автор исправляет недочеты, иногда выясняется известность результата и / или %возможность его развития...
%Обычно в результате мы даем рекомендацию, а также помогаем автору выложить представленную в Архив версию на популярный ресурс по данной теме ... --- на случай задержки или отклонения модераторами Архива.
%Некоторые работы выложены (и могут быть выложены) в Архив независимо от программного комитета --- на их рекомендателях и лежит ответственность за рекомендацию.}
В \cite[(1) в п. 3.1 на стр. 240]{Be14} не приводится (известная автору) ссылка \cite[\S2]{P}, где приводится сходный текст о требовании выкладывания в архив --- 
\emph{c необходимыми предостережениями и только для номинации научно-исследовательских работ}. 
%Видимо, автор понимал, что 
Часто авторы не дают (известных им) необходимых ссылок на подтверждение своей точки зрения, поскольку понимают, что ссылки эту точку зрения опровергают.  
 
%Приведенную ниже цитату  стоит сравнить с еще одним примером необоснованной критики ММКШ, см. 
%Эта критика завуалирована (ММКШ явно не упоминается) и не послана руководству ММКШ - 

\begin{blud}\label{b:number}
\emph{Важнейшим} критерием успешности конференций является \emph{количество} работ на этой конференции. 

Например, в работе \cite{Be14} о критериях успешности конференций (в форме рекомендаций) написано только в  \cite[п. 3.1 и 3.2]{Be14}  
(в \cite[внизу стр. 234 и далее]{Be14} даются рекомендации не по самим конференциям, а по другой работе). 
Многие рекомендации и требования неясны и противоречивы (ибо приведены без конкретики), а одно из них основано на заведомо ложной информации. 
См. подробнее заблуждениe \ref{b:arxiv} и текст после него, а также заблуждениe \ref{r:recom}. 
%Требований к \emph{правильности} работ нет совсем.   
Поэтому на первый план выходит критерий успешности, четко сформулированный в \cite[п. 3.1]{Be14}: 

<<\emph{9....Разумное число докладов. Оптимальное --- примерно 20-25, а если их меньше 10, то это, скорее, неудача.}>>
\end{blud}

Это заблуждение (без его явной публичной формулировки) используется при принятии решений. 
Например, отказ в спонсорской поддержке может быть мотивирован \emph{только} количеством участников. 
Неразумность этого заблуждения очевидна: лучше одному человеку принести пользу, чем десяти --- вред.  
Неразумность ясна еще и из возможности \emph{имитации}: количество участников в моей конференции я могу увеличить,  

$\bullet$ не замечая серьезные недочеты в работах (тогда не получится публиковать полные тексты награждаемых работ, но это общепринято); 

$\bullet$ подавляя распространение информации о других основах исследовательской работы со школьниками, см. \cite[\S1]{ZS26}.  
%\ref{s:intr}.   

%Мебиус!!!
 
Заблуждение \ref{b:number} тем более естественно, что оно близко к разумному соображению: \emph{(достаточно большое) количество школьников действительно должно быть 
важным второстепенным критерием успешности, рассматриваемым с учетом обстоятельств и более важных критериев}.  
Ср. \cite[\S2 Введение: о критериях успешности]{KS}. 
\emph{Важнейшим} же критерием успешности для нас является принцип <<не навреди>>, т.е. соответствие между уровнем содержательности / надежности   награждаемых работ и формулировками награждения, см. конкретику в \S\ref{why}-\S\ref{s:howto}. 
%\cite{V2}.  
Мы тем более благодарны спонсорам и коллегам за материальную и моральную поддержку этого принципа. 
 
\begin{blud}\label{r:recom}
Организаторы конференций призывают к мягкости, но 
%выставляем жесткие требования 
не замечают, что поступают жестко.  
Организаторы конференций выставляют жесткие требования и пишут о высоком уровне конференций, но не замечают, что не выполняют этих 
%(или других, более важных) 
требований. 

Например, в первые годы ММКШ 
%\cite{V} 
иногда общение с авторами и руководителями происходило недостаточно доброжелательно. 
Извинения принесены и приняты.  
В качестве одной из мер по преодолению этого рецензии стали выкладываться на сайт. 

%Вот другой пример. 
В \cite[п. 3.1, (2,3,6,7,8)]{Be14} на разные лады  рекомендуется \emph{мягкость} в отношении докладов, а в \cite[п. 3.2]{Be14} приводятся противоречащие 
%жестко сформулированные 
\emph{требования}.
%Никакой конкретики не приведено. 
%Поэтому рекомендации противоречат требованиям. 
Вот конкретные примеры. 

В \cite[п. 3.1]{Be14} написано: 

<<\emph{3. ... Предварительный отсев докладов должен быть максимально мягким, надо твердо отсекать лишь <<уфологию>> и совсем слабые работы. 
Следует доверять известным специалистам и принимать работы по их рекомендации.}>>

<<\emph{7. ... Следует напоминать жюри, что участники могут быть чрезвычайно чувствительны к оценкам их
выступлений, и призывать к осторожности.}>>

<<\emph{8. Оценка учащегося должна определяться позитивными аспектами:
творческими достижениями, красивыми идеями (в том числе методическими), оригинальностью изложения...}>> 

А в \cite[п. 3.2, сноска 14]{Be14} написано: 

<<\emph{Автор неоднократно сталкивался с ситуацией, когда школьникам и руководителю указывали на математическую ошибку. 
Однако она не исправлялась, а докладчик на следующей конференции повторял доклад слово в слово. 
По мнению автора, это основание для дисквалификации. 
То же относится и к злостному неуказанию источников.}>>

Мы полностью согласны с предыдущей фразой сноски: <<\emph{Нарушением научной этики является отказ в исправлении [ошибок].}>> 
Однако доведение до <<\emph{дисквалификации}>> нам представляется излишне жестким. 

Во-первых, учащийся мог не понять замечания. 
Или не согласиться с ним, и не иметь возможности обсудить несогласие до конференции, а не на ней (или не знать, что заранее лучше). 
Поэтому вместо <<\emph{дисквалификации}>> более мягко не принимать работу на конференцию, пока все важные ошибки не исправлены 
(преодолевая заблуждение \ref{b:number}), при этом давая учащемуся достаточные возможности для их исправления. см. \S\ref{s:rest}). 
%\cite[\S6]{V2}.  
Вместо оценок типа <<злостное>> мы употребляем более спокойные выражения (<<повторное>>).  

Во-вторых, раз работа докладывается на конференции, то она прошла <<\emph{предварительный отсев докладов}>>, к <<\emph{максимальной мягкости}>> 
которого призывает автор (возможно, <<\emph{работа принята по рекомендации известного специалиста}>>). 
Стоит ли наказывать учащегося <<\emph{дисквалификацией}>> за непродуманный автором баланс между желаниями мякости и требовательности? 
Ср.  \cite[\S3.1, сноска 13]{Be14}. 

В третьих, работа с ошибкой вполне  может содержать <<\emph{позитивные аспекты: творческие достижения, красивые идеи, оригинальность изложения}>>, 
Автор призывает именно ими <<\emph{определять оценку [работы] учащегося}>>, что противоречит <<\emph{дисквалификации}>>. 

В четвертых, см. пример в конце \S\ref{s:nomi}. 
%\cite[пример в конце \S5]{V2}. 
  
В \cite[п. 3.1]{Be14} написано: 

<<\emph{2. ... Доброжелательность к учащимся. 
Жюри не должно самоутверждаться за их счет, в частности, проявлять снобизм. 
Необходима доступность в общении, не следует строить из себя небожителей.}>>

А в \cite[п. 3.2]{Be14} написано:

<<\emph{5. Не следует публично решать квадратные уравнения! 
Доклад — не урок, где автор демонстрирует свои элементарные технические умения с подробным рассказом у доски.}>>

Здесь выражение <<квадратные уравнения>> используется не буквально, а как негативный образ для чего-то тривиального.  
%(с точки зрения члена жюри).
\footnote{Действительно, доклады о школьных квадратных уравнениях нам неизвестны (и в \cite{Be14} ссылки на них не приводится). При этом работы о квадратных уравнениях в числовых системах, или о свойствах квадратного трехчлена --- например, динамических --- могут быть весьма нетривиальны.} 
Такие выражения нам представляются недостаточно доброжелательными (снобистскими). 
То, что  тривиально, элементарно, технично с точки зрения члена жюри, может быть серьезным достижением для школьника. 
Даже если это не так, мы 
%считаем желательным 
используем более спокойные выражения (<<молодец, что сделал эту работу, но она не слишком выходит за пределы школьной программы>>, 
<<молодец, что привел подробности в тексте, но в докладе их лучше пропустить>>).

%!!! рец на ZS25

Интересно, что требований к \emph{компетентности рецензирвания} и к \emph{математической грамотности работ} в \cite[п. 3.1 и 3.2]{Be14} нет совсем.   

%правильности 

Рекомендации \emph{мягкости} без явно указанных (и проиллюстрированых примерами) их ограничений \emph{требованиями} (и наоборот) аналогичны 
%являются бессодержательным и  
противоречивому популистскому призыву <<за снижение налогов и повышение качества бесплатной медицины>>
%<<за все хорошее против всего плохого>>. 
Как конкретнее мы \emph{сочетаем} рекомендации по \emph{мягкости} и  \emph{требовательности}, написано в \S\ref{why}-\S\ref{s:howto}. 
%\cite{V}. 
\end{blud}

\emph{О других заблуждениях (в т.ч. транслированных в \cite[название+, 3.1(4), сноски 12, 13]{Be14}) мы при необходимости напишем здесь позже.} 

%Многое другое в этой уже опубликованной в МП статье столь же некорректно:  ...
%https://www.mathnet.ru/php/person.phtml?option_lang=rus&personid=177663
%https://www.mathnet.ru/rus/person99957

\end{document}